\documentstyle{amsart}

\newcommand{\HH}{\operatorname{H}}
\newcommand{\GL}{\operatorname{GL}}
\newcommand{\SL}{\operatorname{SL}}
\newcommand{\SO}{\operatorname{SO}}
\newcommand{\Un}{\operatorname{U}}
\newcommand{\End}{\operatorname{End}}
\newcommand{\Hom}{\operatorname{Hom}}
\newcommand{\Map}{\operatorname{Map}}
\newcommand{\Lie}{\operatorname{Lie}}
\newcommand{\Trace}{\operatorname{Trace}}
\newcommand{\Spin}{\operatorname{Spin}}
\newcommand{\ext}[1]{\overset{#1}{\wedge}}
\newcommand{\by}[1]{@>#1>>}
\newcommand{\tensor}{\otimes}
\newcommand{\contract}{\lrcorner}
\newcommand{\into}{\hookrightarrow}
\newcommand{\bbC}{{\mathbb C}}
\newcommand{\bbR}{{\mathbb R}}
\newcommand{\bbZ}{{\mathbb Z}}
\newcommand{\cN}{{\mathcal N}}
\newcommand{\cC}{{\mathcal C}}
\newcommand{\cA}{{\mathcal A}}
\newcommand{\cW}{{\mathcal W}}
\newcommand{\cG}{{\mathcal G}}
\newcommand{\cM}{{\mathcal M}}
\newcommand{\triv}{\boldsymbol{1}}
\newcommand{\half}{\frac{1}{2}}
\newcommand{\norm}[1]{\mid\!\mid\!#1\!\mid\!\mid}

\title{Seiberg-Witten invariants--an expository account}
\author{Kapil H. Paranjape}
\address{Institute of Mathematical Sciences\\Chennai 600 113,India}
\email{kapil@@math.imsc.ernet.in}
\author{V. Pati}
\address{Stat-Math Unit\\ Indian Statistical Institute\\Bangalore 560 059, India}
\email{pati@@isibang.ac.in}
\begin{document}
\maketitle
We recall some constructions of spin groups in low dimensions.
\section{Spin Groups}
\subsection{Dimension 3}

	Let $W$ be a vector space of dimension 2. Consider the
representation of $\GL(W)$ on $\End^0(W)$, the space of all traceless
endomorphisms of $W$. There is a natural non-degenerate form $\{,\}$
on $\End^0(W)$ given by
 \[	\{f,g\}  = \Trace_W f\circ g	\]
Moreover, we have a sequence of isomorphisms of representations of
$\GL(W)$,
 \[	\ext{3}\End^0(W) = \ext{4}\End(W) = \ext{4}(W^*\tensor W)
	= (\ext{2}W^*)^{\tensor 2} \tensor (\ext{2}W)^{\tensor 2}
	= \triv	\]
 where $\triv$ denotes the trivial representation. Thus we obtain a
natural homomorphism $\GL(W) \to \SO(\End^0(W))$. Over the complex
numbers this identifies $\GL(2)$ with the ``Cspin'' group of $\SO(3)$.
The subgroup $\SL(2)$ is identified with the spin group.

\subsection{Dimension 4}

	Let $W_{+}$ and $W_{-}$ be two vector spaces of dimension 2
and let $\phi:\ext{2}W_{-}\to\ext{2}W_{+}$ be an isomorphism. Then the
vector space $U=\Hom(W_{+},W_{-})$ is isomorphic to its dual via a map
$B:U\to U^*=\Hom(W_{-},W_{+})$ defined by the identity
 \[	\phi(f(w_{+})\wedge w_{-}) = w_{+} \wedge B(f)(w_{-})	\]
 Thus we have a non-degenerate pairing
 \[	(f,g) = \Trace_{W_{+}} B(f)\circ g = \Trace_{W_{-}} g\circ B(f) \]
 which can be seen to be a symmetric form. The group of automorphisms
of the triple $(W_{+},W_{-},\phi)$ is 
\[	S(\GL(W_{+})\times\GL(W_{-})) =
	\{ (g,h) \mid \det(g) = \det(h) \}
\]
We have a sequence of isomorphisms of representations of this group
\[	\ext{4}(U) = \ext{4}(W_{+}^*\tensor W_{-}) 
		= (\ext{2}W_{+}^*)^{\tensor 2} \tensor
			(\ext{2}W_{-})^{\tensor 2}
				\by{\phi\tensor\phi} \triv
\]
Thus we obtain a morphism
\[	S(\GL(W_{+})\times\GL(W_{-})) \to \SO(\Hom(W_{+},W_{-})) \]
Over the complex numbers this identifies the group
$S(\GL(2)\times\GL(2))$ with the ``Cspin''
group of $\SO(4)$ and the subgroup $\SL(2)\times\SL(2)$ is identified
with the spin group of $\SO(4)$.

\subsection{Dimension 6}

Let $U$ be a four dimensional vector space and let
$\psi:\ext{4}U\to\triv$ be a chosen isomorphism so that the group of
automorphisms of the pair $(U,\psi)$ is $\SL(U)$. Consider the pairing
$<,>$ on $\ext{2}U$ given by the composite
 \[	\ext{2}U\tensor\ext{2}U \by{\wedge} \ext{4}U \by{\psi} \triv \]
 This is symmetric and non-degenerate. Moreover, we have a natural
sequence of isomorphisms of representations of $\SL(U)$
 \[ \ext{6}(\ext{2}U) = (\ext{4}U)^{\tensor 6} \by{\psi^{\tensor 6}}
		\triv	\]
Thus we have a representation of $\SL(U)$ in $\SO(\ext{2}U)$. Over the
complex numbers this identifies $\SL(4)$ with the spin group of
$\SO(6)$. 

\subsection{Combination of the above}

Now consider the situation of (1.3) where $U=\Hom(W_{+},W_{-})$. In
this situation $U$ carries a non-degenerate pairing $(,)$ as described
above and hence there is an induced pairing on $\ext{2}U$ which we
also denote by $(,)$. We then have an automorphism $*$ on $\ext{2}U$
defined by the identity $(\alpha,\beta)=<\alpha,*\beta>$. Now the fact
that $\psi=\phi\tensor\phi$ satisfies $(\psi,\psi)=1$ implies that
$*^2=1$. Moreover, one can see that the positive (resp. negative)
eigenspace $\Lambda^{+}$ (resp. $\Lambda^{-}$) of $*$ is of dimension
$3$. Thus the combined representation
 \[ S(\GL(W_{+}) \times \GL(W_{-})) \to \SO(U) \into \SL(U)
		\to \SO(\ext{2}U) \]
 gives a morphism into $\SO(\Lambda^{+})\times\SO(\Lambda^{-})$. Now
we have natural maps $S(\GL(W_{+})\times\GL(W_{-}))\to \GL(W_{\pm})$.
And hence we have representations of $S(\GL(W_{+})\times\GL(W_{-}))$
into $\SO(\End^0(W_{\pm}))$. Consider the homomorphisms of
representations of $S(\GL(W_{+})\times\GL(W_{-}))$
 \[	\ext{2}U \to \End^0(W_{+}) \text{ where } 
		f\wedge g \mapsto B(f)\circ g - B(g) \circ f	\]
and similarly
 \[	\ext{2}U \to \End^0(W_{-}) \text{ where } 
		f\wedge g \mapsto f\circ B(g) - g \circ B(f)	\]
These induce isomorphisms of $\End^0(W_{\pm})$ with $\Lambda^{\pm}$.

\subsection{Compact forms}

Let us fix hermitian structures $h_{\pm}$ on $W_{\pm}$ so that $\phi$ is
an isometry. The group of automorphisms then becomes
$S(\Un(W_{+})\times\Un(W_{-}))$. We define a $\bbC$-anti-linear
automorphism $f\mapsto f^{\dagger}$ defined by the identity
 \[	h_{+}(f^{\dagger}(w),w') = h_{+}(w,B(f)(w'))	\]
One sees that $f^{\dagger\dagger}=f$. Thus we obtain a real vector
space $T$ so that $U= T + \iota T$. Moreover, one sees that the form
$(,)$ restricts to a positive definite form on $T$; hence we obtain a
representation $S(\Un(W_{+})\times\Un(W_{-}))\to \SO(T)$. 
The above discussion then gives us a decomposition of $\ext{2}T$ into
$\Lambda^{\pm}_{\bbR}$. 

We have a $\bbC$-anti-linear endomorphism $f\mapsto f^{\dagger}$ of
$\End^0(W_{\pm})$ given by
 \[	h_{\pm}(f^{\dagger}w,w') = h_{\pm}(w,f(w'))	\]
 One shows that under the isomorphism between $\End^0(W_{\pm})$ and
$\Lambda^{\pm} = \Lambda^{\pm}_{\bbR} + \iota \Lambda^{\pm}_{\bbR}$,
we obtain identifications of $\Lambda^{\pm}_{\bbR}$ with the spaces
$\End^0(W_{\pm})^{ah}$ consisting of $f=-f^{\dagger}$. 

We note that for any pair of elements $\Phi$,$\Psi$ of $W_{+}$ we have an
element $\sigma(\Phi,\Psi)$ of $\End^0(W_{+})$ given by
 \[	w \mapsto i(h_{+}(w,\Psi)\cdot\Phi - \half h_{+}(\Phi,\Psi)\cdot w) \]
When $\Phi=\Psi$ this is an element of $\End^0(W_{+})^{ah}$.
We identify this with an element of $\Lambda^{+}_{\bbR}$.

\subsection{Unitary group case}
	We now further specialise to the case when
$W_{+}=\triv\oplus\det W_{-}$. For ease of notation we use $W$ for
$W_{-}$. In this case, we have a natural sequence of identifications
 \[	\Hom_{\bbC}(W_{+},W_{-})=W\oplus W^* = W\oplus \overline{W}
			= W\tensor_{\bbR}\bbC	\]
Thus we can identify the special orthogonal representation $T$ with
the underlying real vector space of $W$. Now let $\ext{(2,0)}T$ denote
the underlying real vector space to $\ext{2}_{\bbC}W$ and let
$\ext{(1,1)}T$ the real vector space such that
$W\tensor\overline{W}=\ext{(1,1)}T\tensor{\bbR}\bbC$.
We have a natural decomposition
 \[	\ext{2}T = \ext{(2,0)}T \oplus \ext{(1,1)}T	\]
The imaginary part of the hermitian metric on $W$ gives a natural
element $\omega$ of the latter space. One then computes that
 \[ \lambda^{+}_{\bbR} = \ext{(2,0)}T \oplus \bbR\cdot\omega
 \text{ and } \lambda^{-}_{\bbR} = \omega^{\perp} \cap \ext{(1,1)}T \]
Moreover, under the identification between
$\Lambda^{+}_{\bbR}$ and $\End(W_{+})^{ah}$ we obtain identifications
 \[ \ext{(2,0)}T = \Hom_{\bbC}(\triv,\det W) = \det W \text{ and }
	\bbR = \bbR\cdot\omega = \bbR i\cdot\triv_{\det W}		\]

\section{Spin structures on four manifolds}

        Let $X$ be a compact oriented four manifold. For any metric
$g$ on $X$ we have the principal $\SO(4)$ bundle $P$ on $X$ which
consists of oriented orthonormal frames. This corresponds to a class
$[P]$ in $\HH^1(X,\SO(4))$. Using the exact sequence
 \[        1 \to \Un(1) \to \Spin_c(4) \to \SO(4) \to 1        \]
we see that we have an exact sequence
 \[        \HH^1(X,\Un(1)) \to \HH^1(X,\Spin_c(4)) \to
                \HH^1(X,\SO(4)) \to \HH^2(X,\Un(1))        \]
we see that the obstruction to giving a reduction of structure group
from $\SO(4)$ to $\Spin_c(4)$ is given by a class in
$\HH^2(X,\Un(1))$. Moreover, from the exact sequence
 \[        1 \to \bbZ/2\bbZ \to \Spin(4) \to \SO(4) \to 1        \]
we see that the obstruction to giving a spin structure lies in
$\HH^2(X,\bbZ/2\bbZ)$. Under the natural inclusion of $\bbZ/2\bbZ$ in
$\Un(1)$, the obstruction for spin maps to the obstruction for Cspin.
 In fact consider the diagram
\[
\begin{array}{rcccccl}
        & 1                &        & 1 \\
        & \downarrow        &        & \downarrow \\
        1 \to & \bbZ/2\bbZ& \to & \Spin(4) & \to & \SO(4) & \to 1 \\
        & \downarrow        &        & \downarrow &        & \mid\mid\\
        1 \to & \Un(1) & \to & \Spin_c(4) & \to & \SO(4) & \to 1 \\
        & \downarrow        &        & \downarrow \\
        & \Un(1)        & = & \Un(1) \\
        & \downarrow        &        & \downarrow \\
        & 1                &        & 1
\end{array}
\]
 By the associated diagram of cohomologies
\[
\begin{array}{rcccccccl}
        & \downarrow        &        & \downarrow         &        &                &        & \downarrow\\
\to & \HH^1(X,\bbZ/2\bbZ)& \to & \HH^1(X,\Spin(4)) & \to & \HH^1(X,\SO(4)) & \to & \HH^2(X,\bbZ/2\bbZ) \\
        & \downarrow        &        & \downarrow &        & \mid\mid &       & \downarrow \\
\to & \HH^1(X,\Un(1)) & \to & \HH^1(X,\Spin_c(4)) & \to & \HH^1(X,\SO(4)) & \to & \HH^2(X,\Un(1))\\
        & \downarrow        &        & \downarrow  &                &                &        & \downarrow \\
        & \HH^1(X,\Un(1))& = & \HH^1(X,\Un(1)) &        &                &        & \HH^2(X,\Un(1))\\
        & \downarrow        &        & \downarrow 
\end{array}
\]
 we see that the distinct lifts of a given $\SO(4)$ bundle to a
$\Spin_c(4)$ bundle correspond exactly to the different lifts of the
$\Spin(4)$ obstruction class in $\HH^2(X,\bbZ/2\bbZ)$ to a class in
$\HH^1(X,\Un(1))$. We note that the latter is the group of metrised
complex line bundles.

Now we have a natural exact sequence (the exponential sequence) of
sheaves
 \[        0 \to \bbZ \to {\mathcal C}^{\infty} \to \Un(1) \to 1 \]
which gives the natural isomorphisms
$\HH^{i+1}(X,\bbZ)=\HH^i(X,\Un(1))$. Moreover, under these
isomorphisms the exact sequence
 \[ \to \HH^1(X,\Un(1)) \to \HH^1(X,\Un(1)) \to \HH^2(X,\bbZ/2\bbZ)
        \to \HH^2(X,\Un(1)) \to \HH^2(X,\Un(1))                \]
is the same as the exact sequence
 \[ \to \HH^2(X,\bbZ) \to \HH^2(X,\bbZ) \to \HH^2(X,\bbZ/2\bbZ)
        \to \HH^3(X,\bbZ) \to \HH^3(X,\bbZ)                \]
 To summarise, the obstruction to giving a $\Spin_c(4)$ structure is
the image in $\HH^3(X,\bbZ)$ of the obstruction to a $\Spin(4)$ which
lies in $\HH^2(X,\bbZ/2\bbZ)$. If the former is zero then the
different $\Spin_c(4)$ structures correspond to the different lifts of
the $\Spin(4)$ obstruction class to $\HH^2(X,\bbZ)$.
 
In the case when the principal bundle is the one associated with the
metrised tangent bundle as above we have the result that the
obstruction to having a spin structure is given by $w_2(X)$ in
$\HH^2(X,\bbZ/2\bbZ)$; the second Stiefel-Whitney class of $X$. Then
we have Wu's formula which implies that for any
$y$ in $\HH^2(X,\bbZ/2\bbZ)$ we have $w_2(X)\cap y = y\cap y$. Now
consider the image $w$ of $w_2(X)$ in $\HH^3(X,\bbZ)$; this is a
2-torsion class. Let $\HH^2(X,\bbZ)_{\tau}$ denote the group of
torsion elements in $\HH^2(X,\bbZ)$.
There is a natural duality between the 2-torsion in
$\HH^3(X,\bbZ)$ and the group $\HH^2(X,\bbZ)_{\tau}\tensor\bbZ/2\bbZ$;
this duality is given as follows. Let $a\in\HH^2(X,\bbZ)_{\tau}$ be a torsion
class and let $b\in\HH^3(X,\bbZ)$ be a 2-torsion class. Let $b'$ be a
class in $\HH^2(X,\bbZ/2\bbZ)$ whose image is $b$. Let $a'$ be the
image of $a$ in $\HH^2(X,\bbZ/2\bbZ)$, then $<a,b>=(a',b')$. By this
identification we have 
 \[        <a,w>=(a',w_2(X))=\Trace(a'\cap w_2(X)) = \Trace(a'\cap
a')=0\]
for all $a$ in $\HH^2(X,\bbZ)_{\tau}$. But then by the
duality we see that $w$ is 0. Hence, in this case we obtain that
$w_2(X)$ is the reduction modulo 2 of an integral cohomology class; in
other words {\em an oriented compact Riemannian four manifold always
has a $\Spin_c(4)$ structure}.

\section{Monopole equations and their moduli space}
	In this section we describe the monopole moduli spaces and
compute the expected dimension.

\subsection{Connections for Spin structures}
Let $(X,g,c)$ be a compact oriented Riemannian four manifold
with a $\Spin_c$ structure (denoted by $c$). Let $Q$ denote the
corresponding principal $\Spin_c(4)$ bundle over $X$. Then the
principal bundle of oriented orthonormal frames on $X$ is given by
$P=Q/\Un(1)$. We have a natural torsion free connection on this bundle
called the Riemannian connection. The pull-back of this to $Q$ gives
us a 1-form on $Q$ with values in $\Lie(\SO(4)$ which is invariant for
the action of $\Spin_c(4)$. Now consider the principal $\Un(1)$ bundle
$Q/\Spin(4)$ associated with $Q$ which is just the space of all unit
vectors in the line bundle $L=\det(W_{\pm})$. Let $A$ be connection on
this line bundle.  We can pull this back to a form on $Q$. Adding the
above two forms together we obtain a connection on $Q$ which we shall
denote by $\nabla_A$ since the Riemannian connection is unique
whereas $A$ can be varied.

\subsection{Dirac equation}

	Fixing $A$ for the time being we have the differential
operator $\nabla_A:W_{+}\to W_{+}\tensor T^*X$ induced by the
connection as above. One the one hand the Riemannian structure gives
us a natural (flat) identification between $T^*X$ and $TX$ and on the
other we have seen that $TX$ can be thought of as a subspace of
$\Hom_{\bbC}(W_{+},W_{-})$; moreover, this identification is also
invariant under the connection (flat). Thus by contraction we obtain
the composite differential operator of order 1
 \[	D_A = D_{A,+}: W_{+} \to W_{+} \tensor T^*X
		\to W_{+} \tensor TX \to W_{-} \]
 This is called the Dirac operator. As seen earlier we have a natural
identification $\Hom_{\bbC}(W_{+},W_{-})=\Hom_{\bbC}(W_{-},W_{+})$.
Thus we also obtain an operator $D^*_A=-D_{A,-}:W_{-}\to W_{+}$. We have the
identity (see Section 4)
 \[	\int_X h_{-} (-D_{A,+} \Phi,\Psi) 
			= \int_X h_{+} (\Phi, D_{A,-}\Psi) 	\]
so that we see that $D^*_A$ is the adjoint of $D_A$. This
justifies the notation.

The first monopole equation is the Dirac equation $D_A(\Phi)=0$.

\subsection{Second monopole equation}
	Consider the curvature $F_A$ of the connection $A$ on $L$. This
gives a two form with values in the Lie algebra of $\Un(1)$ which is
just $\bbR$. Let $F_A^{+}$ denote the projection into $\Lambda^{+}$.
We have also defined the map $\sigma:W_{+}\tensor\overline{W_{+}}
\to \End(W_{+})^{ah}$. Moreover, we have obtained an
identification of $\Lambda^{+}$ with the space of skew-Hermitian
endomorphisms of $W_{+}$. The second monopole equation is
\[	F_A^+ = \sigma(\Phi,\Phi)	\]

\subsection{Gauge group}
	Let $\cG=\Map(X,\Un(1))$ and consider the action of this group
on the space $\cN$ of pairs $(A,\Phi)$ where $A$ is a connection on
the line bundle $L$ and $\Phi$ a section of $W_+$ given by
 \[	g\cdot(A,\Phi) = (g\cdot A,g\cdot\Phi)
			 = (A - (1/2\pi i)g^{-1}dg, g\Phi)	\]
We see easily that if $(A,\Phi)$ satisfies the monopole equation then
so does $g\cdot(A,\Phi)$. In fact we have
 \[	D_{g\cdot A} (g\cdot \Phi) = g\cdot D_A\Phi	\text{ and }
	F_{g\cdot A} = F_A \text{ and } \sigma(g\Phi,g\Psi) =
\sigma(\Phi,\Psi) \]
 Thus we may consider the `moduli space' of monopoles
 \[	M = M_c = \{ (A,\Phi) \mid D_A\Phi=0 \text{ and }
		F_A^+ = \sigma(\Phi,\Phi) \} / \cG	\]
We will show that for `good' metrics this is a compact orientable
manifold. We shall also find out how it depends on this choice of
metric.

	Now let $\cW_{\pm}$ be the spaces of sections of $W_{\pm}$. 
We have a map
 \[	\nu:\cN \to \cW_{-} \text{ given by } 
		(A,\Phi) \mapsto (D_A\Phi)	\]
 Let $\cM$ denote the inverse image $\nu^{-1}(0)$ and let $\cM^*$
denote the open subset consisting of pairs $(A,\Phi)$ where $\Phi\neq
0$. The differential of the map $(A,\Phi)\to D_A\Phi$ is given by
 \[	(a,\phi) \mapsto D_A\phi + 2\pi i a\circ\Phi	\]
Suppose $\Psi$ is orthogonal to the image. Then we obtain the equations
 \[	D_A^*\Psi=0 \text{ and }  \Phi\tensor\Psi = 0	\]
by orthogonality with the image of vectors of the form $(a,0)$ and
$(0,\phi)$ respectively. Now a solution of an elliptic operator
vanishes on an open set only if it is identically 0. Thus we see that
$\psi=0$; in other words $\nu$ is a submersion when restricted to the
space $\cN^*$ consisting of pairs $(A,\Phi)$ where $\Phi\neq 0$.
Thus $\cM^*$ is a manifold (albeit of infinite dimension).

The group $\cG$ acts freely on $\cM^*$ since a solution of an elliptic
operator cannot vanish on an open set unless it is 0. Consider the space
$\Omega^{2+}$ consisting of 2-forms invariant under $*$ with the trivial
action of $\cG$.  The map $\cM\to\Omega^{2+}$ given by
$F_A^{+}-\sigma(\Phi,\Phi)$ factors thorugh the quotient $\cM/\cG$.  We
thus obtain a `complex' $\cG\to \cM^*\to\Omega^{2+}$. The moduli space
can be thought of as being its `cohomology'.

\subsection{Virtual dimension of the moduli space}
	To compute the dimension of the moduli space we need to
compute the cohomology of the complex of differentials of the complex
$\cG\to\cM^*\to\Omega^{2+}$. The tangent space to $\cG$ at identity
can be identified with $\Omega^0$ the space of fucntions and the
tangent space to $\Omega^{2+}$ can be identified with itself since it
is a vector space. We have an exact sequence 
 \[ 0 \to T\cM^* \to \Omega^1 \oplus \cW_{+} \to \Omega^{2+}\to 0\]
Where we have identified the tangent space of $\cA$ with $\Omega^1$
the space of 1-forms.  Thus the complex of differentials 
 \[ T\cG \to T\cM^* \to \Omega^{2+}	\] 
is quasi-isomorphic to the complex
 \[	\Omega^0 \to \Omega^1 \oplus \cW_{+} \to
		\Omega^{2+} \oplus \cW_{-}	\]
 where the maps are
 \[	h \mapsto (-dh, 2\pi i h\Phi) \text{ and } 
	(a,\psi) \mapsto (d^{+}a - \Im\sigma(\Phi,\phi),
				D_A\phi + 2\pi i a\circ \Phi)	\]
for $h\in\Omega^0$, $a\in\Omega^1$ and $\phi\in\cW_{+}$. Here $d^{+}$
denotes the exterior derivative combined with the projection to
$\Omega^{2+}$ and $\Im\sigma(\Phi,\phi)$ denotes the skew-Hermitian part of
$\sigma(\Phi,\phi)$. This complex is homotopic to the complex where the
first map is $h\mapsto (-dh,0)$ and the second is
$(a,\phi)\mapsto(d^{+}a,D_A\phi)$ since the difference between these
two complexes is given by compact operators. Thus the index of our
complex of differentials is the index of the complex
 \[	\Omega^0 \to \Omega^1 \oplus \cW_{+} \to
		\Omega^{2+} \oplus \cW_{-}	\]
 where the maps are
 \[	h \mapsto (-dh, 0) \text{ and } 
	(a,\psi) \mapsto (d^{+}a , D_A\phi )	\]
 This is a topological invariant for the pair $(X,c)$ by the Atiyah-Singer Index
theorem; we call this the virtual dimension of the moduli space.
In case we can find a point $\delta\in\Omega^{2+}$ which is a
regular value for the map $\cM^*\to\Omega^{2+}$ we see that this index
will be the dimension of
 \[	M_{c,\delta} = \{ (A,\Phi) \mid
		D_A\Phi=0 \text{ and } F_A^{+} = \sigma(\Phi,\Phi)
					+\delta	\}/cG	\]
We call this the perturbed moduli space. We will show that such a
value of $\delta$ exists that $M_{c,\delta}$ is a compact
orientable manifold whose dimension is the virtual dimension.

\section{Differential Calculus}
	We derive various identities among differential operators in
the context of $\Spin_c$ connections.

\subsection{The Adjoint of the Dirac operator}

	We have defined the Dirac operator as the composite
\[	D_{A,+} = D_{A}: W_{+} \by{\nabla_A} TX^*\tensor W_{+}
		\to	W_{-}	\]
 where the latter map is the contraction under the identification of
$TX^*$ with $TX\subset\Hom_{\bbC}(W_{+},W_{-})$. We have similarly the
Dirac operator $D_{A,-}: W_{-} \to W_{+}$ since we have an
identification of $\Hom_{\bbC}(W_{+},W_{-})$ with its dual space
$\Hom_{\bbC}(W_{-},W_{+})$. In terms of an orthonormal frame of
tangent vectors $\{e_i\}$ we obtain a sequence of identities:
 \[ (D_A\Phi,\Psi) = \sum_i (e_i\circ\nabla_{e_i}\Phi,\Psi) \]
 and since $(f\circ\Phi,\Psi) =
(f^{\dagger}\circ\Phi,\Psi)=(\Phi,f\circ\Psi)$ for all $f$ in $TX$,
 \[ (D_A\Phi,\Psi) = \sum_i (\nabla_{e_i}\Phi,e_i\circ\Psi) \]
 Now the fact that $\nabla$ is a metric connection means that
 \[	(\nabla_{e_i}\Phi,e_i\circ\Psi) = e_i(\Phi,e_i\circ\Psi)
		- (\Phi,\nabla_{e_i}(e_i\circ\Psi))	\]
 Let $d\tau$ denote the volume form then for any function $f$
and any vector field $v$ we have,
 \[	v(f)d\tau = d(f (v\contract d\tau)) - f d(v\contract d\tau) \]
 Thus we obtain
{small
 \[ (D_A\Phi,\Psi)d\tau = \sum_i d((\Phi,e_i\circ\Psi)e_i\contract d\tau)
		- \sum_i (\Phi,e_i\circ\Psi)d(e_i\contract d\tau)
		- \sum_i (\Phi,\nabla_{e_i}(e_i\circ\Psi))	\]
}
 For any vector field $v$ we have the identity
 \[	d(v\contract d\tau)  = \sum_j (e_j,\nabla_{e_j}v)d\tau \]
  Moreover, since $(e_j,e_j)=\delta_{i,j}$ is a constant we have
{\small \[	\sum_i (\Phi,e_i\circ\Psi)d(e_i\contract d\tau) =
 	  \sum_{i,j} (\Phi,e_i\circ\Psi) (e_j,\nabla_{e_j}e_i) d\tau
	= - \sum_{i,j} (\Phi,e_i\circ\Psi) (\nabla_{e_j}e_j,e_i) d\tau
\]}
 The other term can be written as follows
 \[ \nabla_{e_i}(e_i\circ\Psi) = (\nabla_{e_i}e_i)\circ\Psi 
		+ e_i\circ\nabla_{e_i}\Psi \]
 and 
 \[	(\nabla_{e_i}e_i)\circ\Psi =
		\sum_j (\nabla_{e_i}e_i,e_j)e_j\circ\Psi	\]
 Combining the above identities we obtain
 \[	(D_A\Phi,\Psi) = \sum_i d((\Phi,e_i\circ\Psi)e_i\contract d\tau)
	  - \sum_{i} (\Phi,e_i\circ(\nabla_{e_i}\Psi))d\tau	\]
 Hence 
 \[	\int_X (-D_{A,+} \Phi,\Psi) = \int_X (\Phi, D_{A,-}\Psi) 	\]
 and $-D_{A,-}$ is the adjoint operator of $D_{A,+}$.

By an entirely similar chain of reasoning we show that the adjoint
$\nabla^*:TX\tensor W_{-}\to W_{+}$ of $\nabla$ on $W_{+}$ is given by 
 \[	\nabla(v\tensor\Phi) = -(\sum_i (e_i,\nabla_{e_i}v)\Phi
		+ \nabla_v\Phi)	\]
In invariant terms, we can describe this as the composite
 \[	TX\tensor W_{+} \by{-\nabla} TX^*\tensor TX\tensor W_{+}
		\by{\Trace\tensor \triv} W_{+}		\]

\subsection{The Weitzenbock formula}

	We now compute the composite $D_A^*D_A\Phi$. As before we
choose a local orthonormal frame $\{e_i\}$ for $X$. We then have
 \[ D_A^*D_A\Phi = \sum_i -D_A(e_i\circ \nabla_{e_i}\Phi)
	= - \sum_{i,j} e_j\circ\nabla_{e_j}(e_i\circ\nabla_{e_i}\Phi)) \]
We expand the summand to obtain
 \[	e_j\circ\nabla_{e_j}e_i\circ\nabla_{e_i}\Phi +
		e_j\circ e_i \nabla_{e_j}\nabla_{e_i}\Phi \]
As above the first term above can be expanded again as
 \[ \sum_k (e_k,\nabla_{e_j}e_i) e_j\circ e_k\circ \nabla_{e_i}\Phi =
	- \sum_k (\nabla_{e_j}e_k,e_i) e_j\circ
			e_k\circ\nabla_{e_i}\Phi 	\]
We obtain the formula
 \[ D_A^*D_A\Phi
	= \sum_{i,j,k}(\nabla_{e_j}e_k,e_i)e_j\circ e_k\circ\nabla_{e_i}\Phi
	- \sum_{i,j} e_j\circ e_i\circ\nabla_{e_j}\nabla_{e_i}\Phi \]
 Now defining $\nabla^2_{V,W}=\nabla_V\nabla_W - \nabla_{\nabla_V W}$,
 \[ D_A^*D_A\Phi = - \sum_{i,j} e_j\circ e_i\circ\nabla^2_{e_j,e_i}\Phi \]
 Similar calculations yield the formula
 \[ \nabla_A^*\nabla_A\Phi = - \sum_{i} \nabla^2_{e_i,e_i}\Phi \]
 Now the difference gives us
 \[ D_A^*D_A\Phi - \nabla_A^*\nabla_A\Phi
	= - \sum_{i\neq j} e_j\circ e_i\circ\nabla^2_{e_j,e_i}\Phi \]
 From the definition of $\nabla^2_{\cdot,\cdot}$ we have
{\small \[ \nabla^2_{V,W} - \nabla^2_{V,W}
	= \nabla_V\nabla_W - \nabla_{\nabla_V W} 
		- \nabla_W\nabla_V - \nabla_{\nabla_W V}
	= \nabla_V\nabla_W - \nabla_W\nabla_V - \nabla_{[V,W]}	\]
}
 using the fact that the connection is torsion free. Since we have an
orthonormal basis we have $e_i\circ e_j = - e_j\circ e_i$ so that we
obtain
 \[ D_A^*D_A\Phi - \nabla_A^*\nabla_A\Phi =
	= - \sum_{i<j} e_j\circ e_i\circ R(e_j,e_i)\Phi	\]
where $R(V,W)=\nabla_V\nabla_W-\nabla_W\nabla_V - \nabla_{[V,W]}$ is
the curvature tensor. 

\subsection{The Curvature tensors}
	The $\Spin_c$ connection has been expressed as a sum of the
Riemannian connection and the $\Un(1)$ connection $A$ on $L$. Thus the
curvature tensor $R$ is also the sum of the Riemann curvature tensor
$S$ and the curvature of $A$. The former can be expressed as
 \[	S(V,W)	= \sum_{k,l} (S(V,W)e_l,e_k) e_k\circ e_l	\]
 Thus we obtain
{\small \[ \sum_{i,j} e_j\circ e_i\circ S(e_j,e_i) =
	\sum_{i,j,k,l} (S(e_j,e_i)e_l,e_k)
		e_j\circ e_i\circ e_k \circ e_l		\]}
 By the orthonormality of $e_i$'s we easily resolve the latter to obtain
$\sum_{i,j} (S(e_j,e_i)e_i,e_j)$ which is the negative of the scalar
curvature $s$. The curvature of $A$ considered as an operator on
$W_{+}$ acts as $2\pi i F_A$.  Thus the final (Weitzenbock) formula
reads

 \[	D_A^* D_A - \nabla_A^*\nabla_A = s - 2\pi i F_A 	\]

\subsection{Extrema}
	Let $x$ be a point of our manifold where $(\Phi,\Phi)$ attains
a maximum. Then for any vector $v$ at $x$ we have
$v((\Phi,\Phi))(x)=0$. Thus consider the following identity (where
$\Re$ denotes the real part)
 \[ \Re(\nabla_A^*\nabla_A\Phi,\Phi)
	= - \sum_i \half(e_i e_i (\Phi,\Phi) -
			(\nabla_{e_i}\Phi,\nabla_{e_i}\Phi))
	  + \sum_{i,j} (e_j,\nabla_{e_i}e_i)\Re(\nabla_{e_i}\Phi,\Phi)
 \]
Since $e_i(\Phi,\Phi)(x)=0$ the last term vanishes at $x$. Moreover,
since $x$ is a local maximum for $(\Phi,\Phi)$ the term $e_i
e_i(\Phi,\Phi)(x)$ is negative.  Thus we see that
$\Re(\nabla_A^*\nabla_A\Phi,\Phi)$ is positive at $x$.
\section{The Seiberg-Witten invariants}
	In this section we construct the Seiberg-Witten invariants.
First of all we fix a four manifold $X$, a Riemannian metric $g$ and a
$\Spin_c$ structure $c$. At the end of the section we will discuss the
independence of the invariants on the metric considered.

\subsection{Statement of the basic construction}

	Let $M_{\delta}$ denote the fibre of $\cM/\cG\to\Omega^{2+}$
over the point $\delta$. We wish to show that there is a $\delta$ such
that this is a compact manifold. To show this we need to show
\begin{enumerate}
\item There are regular values for $\cM^*/\cG\to\Omega^{2+}$.
\item There are regular values as above such that the fibre of
$\cM/\cG\to\Omega^{2+}$ is contained in $\cM^*$. 
\item The map $\cM/\cG\to\Omega^{2+}$ is proper.
\end{enumerate}

\subsection{Properness}

Let $\delta_i$ be a convergent sequence of elements in $\Omega^{2+}$.
This means that the sequence converges in the $L^2_k$ Sobolev norm for
every $k$. Let $(A_i,\Phi_i)$ be such that $D_{A_i}\Phi_i=0$ and
$F_{A_i}^+ - \sigma(\Phi_i,\Phi_i) = \delta_i$. To show properness we
need to find a convergent subsequence of $(A_i,\Phi_i)$; for which it
is enough to show that this sequence is bounded in the $L^2_k$ Sobolev
norm for every $k$.

Let $B$ be a fixed smooth connection on $L$; we express $A_i=B+a_i$
where $a_i$ are 1-forms. Consider the function $h_i = G*d*a_i$ and let
$g^0_i=\exp(2\pi i h_i)$. Then $g^0_i\cdot A_i = B + a_i - dh_i$ and
we obtain ${}*d*(a_i - dh_i)=0$. Now we can choose $g_i$ so that the
Hramonic part of $a_i$ lies in the fundamental domain for
$\HH^1(X,\bbZ)$ in $\HH^1(X,\bbR)$. Thus
upto gauge invariance we can replace $A_i$ by
another so that ${}*d*a_i=0$ and the harmonic part $\alpha_i$ of $a_i$
is bounded. Let $b_i=a_i- \alpha_i$. The second monopole equation
becomes
 \[	d^{+} b_i = \sigma(\Phi_i,\Phi_i)-F_B^{+} + \delta_i \]
So that an $L^2_k$ bound on $\Phi_i$ will give an $L^2_k$ bound on
$d^{+} b_i$. But now $b_i = G *d* d^{+} b_i$ by the above construction
of $b_i$; here $G$ is the Green's operator. Thus we obtain a
bound on the $L^2_{k+1}$ norm of $b_i$ since $G$ is 2-smoothing.

Let us write $\Phi_i=\Psi_i+\phi_i$ where $D_B\Psi_i=0$ and $\phi_i$
is orthogonal to the space of solutions of $D_B$. Hence
$\phi=GD_B^*D_B\Phi$ and the first monopole equation becomes 
 \[ D_B\Phi_i = - (b_i + \alpha_i)\circ\Phi_i	\]
so that an $L^2_k$ bound on $\Phi_i$ and $b_i$ gives us an $L_2^{k+1}$
bound on $\phi_i$. We also need to find a way to uniformly bound
$\Psi_i$. We do this by finding a uniform bound for $\Phi_i$.

Let $x_i$ be a point where $(\Phi_i,\Phi_i)$ attains a
supremum. Applying the Weitzenbock formula we see that at $x_i$ we
have 
 \[ 0 \geq \Re(\nabla_{A_i}^*\nabla_{A_i}\Phi_i,\Phi_i)(x_i) =
	-\Re(s\Phi_i,\Phi_i)(x_i) + 2\pi\Im(F_{A_i}\Phi_i,\Phi_i)(x_i) \]
Note that $F_{A_i}$ is a skew-Hermitian endomorphism of $W_{+}$ and
thus 
 \[	\Im(F_{A_i}\Phi_i,\Phi_i) = (F_{A_i}\Phi_i,\Phi_i)	\]
The second monopole equation gives us 
 \[ F_{A_i}\Phi_i = \sigma(\Phi_i,\Phi_i)\Phi_i + \delta_i\Phi_i \]
and the expression for $\sigma$ gives us
 \[ \sigma(\Phi_i,\Phi_i)\Phi_i = {\frac{i}{2}}(\Phi_i,\Phi_i)\Phi \]
Combining the above we obtain 
 \[	(\Phi_i,\Phi_i)(x_i) \leq \max\{0,-s+\norm{\delta_i}\}	\]
Thus we uniformly bound $\Phi_i$ in the $C^0$-norm. This gives us
uniform bounds for $\Psi_i$ and $\phi_i$ in the $C^0$-norm. Now
$\Psi_i$ are solutions of the Dirac equation $D_B\Psi_i=0$. Thus the
set of $C^0$-bounded solutions is a compact set; in particular, we
obtain $L^2_k$ bounds on $\Psi_i$ for all $i$.

The above arguments applied inductively gives the required result.
We note that the above arguments also prove that the solutions of the
monopole equations are smooth since any solution which is bounded in
$L^2_k$ norm for some $k$ is actually bounded in all $L^2_{k'}$ norms
as above.

\subsection{Regular values} 
Now consider the compact space $M_0$ of solutions of the unperturbed
monopole equations. For each point $(A,\Phi)$ of $M_0$ we have a
neighbourhood $U$ in $\cN$ of $(A,\Phi)$ and a finite dimensional
linear space $H\subset\Omega^{2+}$ such that the composite
 \[	U \to \cN \to \cW_{-}\times\Omega^{2+} \to H^{\perp}	\]
is a submersion. By compactness we can find a common $H$ and a
saturated (for $\cG$) open set $U$ in $\cN$ containing the inverse
image of $M_0$ such that the above composite is a submersion. Since
the derivative is a Fredholm map, the fibre over $0$ is a finite
dimensional manifold $N$. We now consider the map of finite
dimensional manifolds $N\to H$. By Sard's theorem we have a dense
subset of $H$ which consists of regular values.

Now assume that $b_2^{+}$ which is the codimension in $\Omega^{2+}$ of
the $\delta$'s of the form $F_B^{+} + d^{+}b$ is greater than zero.
Then the collection of those $\delta$ for which the fibre is contained
in $\cM^*$ is a non-empty open set. If $b_2^{+}>1$ then this open set
is even path-wise connected. Thus in this situation the cobordism
class of the fibre is independent of the regular value chosen.

\subsection{Dependence on the metric}
	Let $\cC$ denote the space of all metrics $g$ on $X$ under
which the fixed volume form $d\tau$ has norm one. We have a natural
map $\cC\to G$ where $G$ denotes the Grassmannian of rank $b_2^+$
quotients of $\HH^2(X,\bbR)$. The corresponding tangent level map is
\[ \Hom(\Lambda^-,\Lambda^+) = T\cC \to TG= \Hom(\HH^{2-}_g,\HH^{2+}_g) \]
Where a map $f:\Lambda^-\to\Lambda^+$ goes to its harmonic projection.

For any class $c=c_1(L)$ in $\HH^2(X,\bbR)$ let $S_c$ denote the
subvariety of $G$ where the class $c$ goes to zero in $\HH^{2+}$. At a
point of $S_c$ the tangent space to $S_c$ is given by the kernel of
the evaluation map $g\mapsto g(c)$. Consider the composite map
 \[ \Hom(\Lambda^{-},\Lambda^{+})  \to \Hom(\HH^{2-}_g,\HH^{2+}_g)
		\to \HH^{2+}_g	\]
 If we show that this map is surjective, then the space of all metrics
under which the class $c$ becomes $*$-anti-invariant will be of
codimension $b_2^{+}$. The argument of the previous section will apply
to show that the Seiberg-Witten invariant is independent of the metric
when $b_2^{+}>1$.

To show that the above map is surjective suppose that $d$ is
perpendicular to the image. We will then obtain that $c\tensor d$ is
identically zero. But now if $c\neq 0$ then it is represented by a
harmonic form which cannot vanish on an open set. Thus $d$ must vanish
on an open set. But we represent $d$ by a harmonic form too. Thus
$d=0$ as required.

\section{The case of K\"ahler manifolds}
	We now specialise to the case of K\"ahler surfaces.
\subsection{Spin structures}
	For any four manifold with almost complex structure and
(hermitian) metric we have a natural $\Spin_c$ structure given by
taking $W^0_{+}=\ext{2}_{\bbC}TX\oplus \triv$ and $W^0_{-}=TX$.
The inclusion of $TX$ in $\Hom_{\bbC}(W^0_{+},W^0_{-})$ is the natural one
as discussed at the end of section 1. Thus any $\Spin_c$ structure on
$X$ is given by $W_{+}=M\tensor_{\bbC}\ext{2}_{\bbC}TX\oplus M$ and
$W_{-}=TX\tensor_{\bbC}M$. For ease of notation we adopt the standard
convention $\ext{2}TX^*=K_X$.

\subsection{Spin$_c$ connections}
	Any $\Un(2)$ connection on $TX$ gives a connection on all
associated bundles. In particular we obtain connections on
$W^0_{\pm}$. However, in order that these be $\Spin_c$ connections it
is necessary that the induced connection on $TX$ be the Riemannian
(torsion-free) connection. This can only happen if the (almost)
complex structure is parallel with respect to the Riemannian
connection; thus in this case the manifold must be K\"ahler.

	To give a connection in the general $\Spin_c$ structure we
need in addition to give a $\Un(1)$ connection on $M$.

\subsection{The First monopole equation}
	Consider a $\Spin_c$ connection as above. We then obtain a
Dirac operator on $M\oplus M\tensor K_X^{-1}$. By the above discussion
we note that the restriction of this to $M$ is the composite
 \[ M \to M\tensor_{\bbR} TX^* =
	M\tensor_{\bbC} TX^* \oplus M\tensor_{\bbC} \overline{TX^*}
		\to M \tensor_{\bbC} TX	\]
Here we have used the identification of $TX^*$ with $\overline{TX}$
given by the hermitian structure. The first map in the above composite
is the $\Un(1)$ connection on $M$. Thus we see that the restriction of
the Dirac operator to $M$ is $\nabla^{(0,1)}$. We similarly show that
the restriction of the Dirac operator to $M\tensor K_X^{-1}$ is also
$\nabla^{(0,1)}$ for the induced $\Un(1)$ connection on this line
bundle. 

\subsection{The Second Monopole equation}
	Following Section 1 we compute that the $(2,0)$ part of
$\sigma(\Phi,\Phi)$ for $\Phi=(\alpha,\beta)$ is
$\overline{\alpha}\beta$ and the $(1,1)$ part of is
$\half(\norm{\beta}^2 - \norm{\alpha}^2)$. Thus the second monopole
equation becomes
 \[ (F_A^{+})^{(2,0)} = {\overline{\alpha}}\beta \text{ and }
	(F_A^{+})^{(1,1)} = \half(\norm{\beta}^2-\norm{\alpha}^2)\omega
\]

\subsection{The Weitzenbock formula}
	We next apply the Weitzenbock formula for any pair $(A,\Phi)$
 \[ D_A^*D_A\Phi = \nabla_A^*\nabla_A\phi +
		s\Phi -2\pi i F_A^{+}\Phi	 \]
to obtain an equality of global inner products
 \[ (D_A\Phi,D_A\Phi)_X = (\nabla_A\Phi,\nabla_A\Phi)_X +
		(s\Phi,\Phi)_X +2\pi \Im(F_A^{+}\Phi,\Phi)_X	\]
 On the other hand we compute the global norm of $F_A
-\sigma(\Phi,\Phi)$ as follows
 \[ \norm{F_A^{+}-\sigma(\Phi,\Phi)}^2_X = \norm{F_A^{+}}^2_X +
	\norm{\sigma(\Phi,\Phi)}^2_X -
		2 \Re(F_A^{+},\sigma(\Phi,\Phi))_X		\]
The last term is computed by the integral of the function
 \[ \Re \Trace_{W_{+}} (F_A^{+}\circ\sigma(\Phi,\Phi))
	= -\Im (F_A^{+}\Phi,\Phi) + \half\Trace(F_A^{+})\norm{\Phi}^2 \] 
 Now $\Trace(F_A^{+})$ is 0. Thus adding the above two identities we
obtain
 \[ \norm{D_A\Phi}^2_X + \norm{F_A^{+}-\sigma(\Phi,\Phi)}^2_X =
	\norm{\nabla_A\Phi}^2_X + (s\Phi,\Phi)_X +
	2\pi(\norm{F_A^{+}}^2_X + \norm{\sigma(\Phi,\Phi)}^2_X) \]
We note that the right hand side is equal to
 \[ \norm{\nabla_A\alpha}^2_X + \norm{\nabla_A\beta}^2_X +
	(s\alpha,\alpha)_X + (s\beta,\beta)_X +
	2\pi\norm{F_A^{+}}^2_X +
	2\pi (\norm{\alpha}^2_X + \norm{\beta}^2_X)^2		\]
which is invariant under a change of sign for $\alpha$ or $\beta$.

Now suppose that $(A,\Phi)$ solve the monopole equations and consider
the pair $(A,\Phi_1)$ where $\Phi_1=(\alpha,-\beta)$. By the above
discussion we see that $(A,\Phi_1)$ is also a solution for the
monopole equations. But then we must have
\[ (F_A)^{(2,0)}=\overline{\alpha}\beta=
			-\overline{\alpha}\beta = 0\]
Thus we obtain the fact that $F_A$ is a holomorphic connection on
$M^{\tensor 2}\tensor K_X$. Moreover, by ellipticity of the Dirac
operator (and its components) we must have that either $\alpha$ or
$\beta$ is zero according as $(F_A^{+})^{(1,1)}$ is a positive or
negative multiple of $\omega$. By the first monopole equation it then
follows that $\alpha$ and $\beta$ are holomorphic sections of the
corresponding line bundles.

\end{document}